\documentclass[11pt]{article}

\usepackage[english]{babel}
\usepackage[latin1]{inputenc}
\usepackage{amsmath,amssymb,enumerate}
\usepackage{ulem}
\usepackage{psfig,labelfig}
\usepackage{graphicx}
\usepackage[dvips]{epsfig}
\usepackage[T1]{fontenc}

\newtheorem{thm}{Theorem}[section]

\newtheorem{lem}[thm]{Lemma}

\newcommand{\C}{{\mathbb C}}
\newcommand{\E}{{\mathbb E}}
\newcommand{\Hy}{{\mathbb H}}
\newcommand{\Sp}{{\mathbb S}}
\newcommand{\R}{{\mathbb R}}

\DeclareMathOperator{\arcsinh}{arcsinh}

\DeclareMathOperator{\iso}{Isom}

\DeclareMathOperator{\diam}{diam}
\DeclareMathOperator{\capa}{cap}
\DeclareMathOperator{\caps}{cap_{\partial s}}
\DeclareMathOperator{\area}{area}

\DeclareMathOperator{\lip}{Lip}

\title{Inequalities for the capacity of non-contractible annuli on cylinders of constant and variable negative curvature}

\author{Muetzel, Bjoern\thanks{E-mail address : bjorn.mutzel@gmail.com}\\
 \\
\small Department of Mathematics, Ecole Polytechnique F\'{e}d\'{e}rale de Lausanne, Station 8, \\
\small CH-1015 Lausanne, Switzerland \\[-0.8ex]
\\
\small Mathematics Subject Classifications (2010):  32Q05, 32U20, 35A15 and 49Q05.}

\begin{document}
\maketitle

\begin{abstract}
Using a new method we give elementary estimates for the capacity of non-contractible annuli on cylinders and provide examples, where these inequalities are sharp. Here the lower bound depends only on the area of the annulus. In the case of constant curvature this lower bound is obtained with the help of a symmetrization process that results in an annulus of minimal capacity. In the case of variable negative curvature we obtain the lower bound by constructing a comparison annulus with the same area but lower capacity on a cylinder of constant curvature. The methods developed here have been applied to estimated the energy of harmonic forms on Riemann surfaces in \cite{mu}.\\

\end{abstract}

\section{Introduction}

\begin{figure}[h!]
\SetLabels
\L(.08*.80) $E_1$\\
\L(.27*.33) $\gamma$\\
\L(.23*.33) $p$\\
\L(.23*.60) $\delta_p$\\
\L(.30*.60) $S$\\
\L(.54*.80) $C$\\
\L(.51*.50) $A$\\
\L(.59*.40) $\gamma$\\
\L(.82*.80) $C$\\
\L(.78*.70) $B$\\
\L(.87*.40) $\gamma$\\
\endSetLabels
\AffixLabels{%
\centerline{%
\includegraphics[height=6cm,width=14cm]{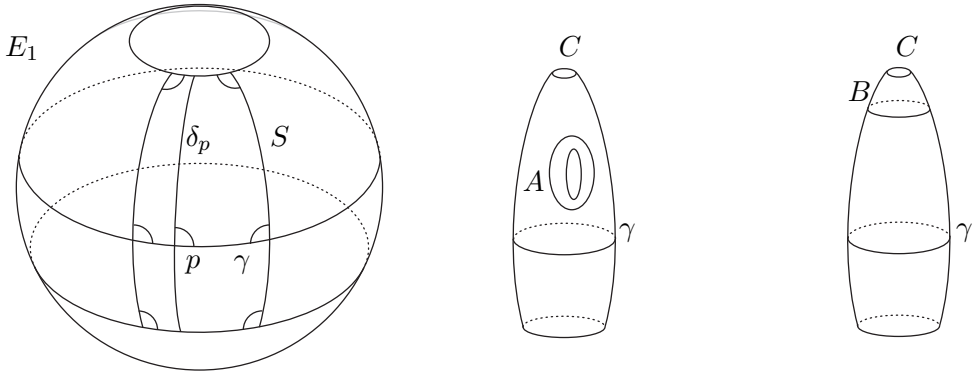}}}
\caption{A strip $S$ in $E_1$ with baseline $\gamma$, the corresponding cylinder $C$ with a contractible annulus $A$ and a non-contractible annulus $B$ that has minimal capacity among all annuli with area $\area(B)$.}
\label{fig:poscurv2}
\end{figure}

\begin{figure}[h!]
\SetLabels
\L(.08*.90) $E_{-1}$\\
\L(.31*.70) $S$\\
\L(.24*.49) $\gamma$\\
\L(.54*.90) $C$\\
\L(.49*.37) $A$\\
\L(.58*.50) $\gamma$\\
\L(.82*.90) $C$\\
\L(.78*.44) $B$\\
\L(.86*.50) $\gamma$\\
\endSetLabels
\AffixLabels{%
\centerline{%
\includegraphics[height=6cm,width=14cm]{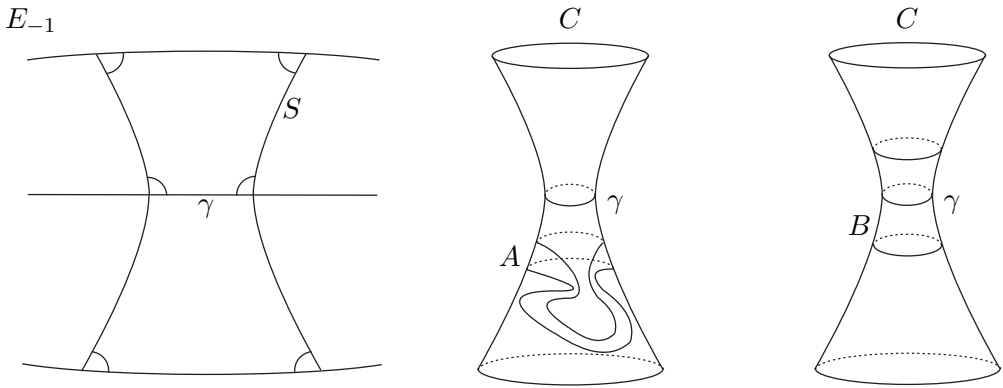}}}
\caption{A strip $S$ with baseline $\gamma$ in a surface $E_{-1}$ with constant curvature $-1$, the corresponding cylinder $C$ with a non-contractible annulus $A$ and the annulus $B$ that has minimal capacity among all annuli with fixed area $\area(A)=\area(B)$.}
\label{fig:negcurv2}
\end{figure}

We will first give the definition of a cylinder. Let $E_K$ be a simply connected surface of constant (sectional) curvature or variable negative curvature, where $K$ denotes the supremum of the curvature in $E_K$.\\
Let $\gamma  \subset E_K$ be a geodesic arc of length $l(\gamma) < \diam(E_K)$. For each point $p \in \gamma$ there exists a geodesic arc $\delta_p$ that is perpendicular to $\gamma$ and that passes through $p$.
We denote by a \textit{strip} $S$ of constant length the set that is formed by the union of such geodesic arcs, i.e.
\[
S= \bigcup \limits_{p \in \gamma} \left\{ \delta_p \right\}.
\]
Here all endpoints of the arcs $\delta_p$ on one side of $S$ shall have the distance $a$ to $\gamma$ and the endpoints on the other side shall have distance $b$ to $\gamma$, where $b \leq a < \frac{\diam(E_K)}{2}$ (see \textit{Fig.~\ref{fig:poscurv2}} and \textit{Fig.~\ref{fig:negcurv2}}). This condition ensures that the different $\delta_p$ do not intersect.\\
A \textit{cylinder with baseline $\gamma$} or shortly \textit{cylinder} $C$ is
\[
C= S \mod M,
\]
such that the sides of $S$ containing the endpoints of  $\gamma$ are identified by an isometry $M \in \iso^+(E_K)$.
Here we assume that such a $M$ exists for $S$, which is always true in the case of constant curvature.

There exist two types of annuli on $C$, contractible and non-contractible annuli (see \textit{Fig.~\ref{fig:poscurv2}, Fig.~\ref{fig:negcurv2}} and section 2). We will consider only non-contractible annuli which we call shortly \textit{annuli} in the following sections. \\
\\
We obtain the following two theorems concerning the lower bound of the capacity of annuli on cylinders:
\begin{thm}
Let $C$ be a cylinder of variable negative curvature smaller than $K=-k^2$. Let $\gamma$ be its baseline of length $l(\gamma)=l$. Let $A \subset C$ be a non-contractible annulus of finite area $\area(A)$. Then for $W'=\arcsinh(\frac{k\area(A)}{2 \cdot l})$ and $h_1(x)=2\arctan(\exp(x))$, we have
\[
     \capa(A) \geq \frac{k \cdot l}{h_1(W')-h_1(-W')}.
\]
\label{thm:main1}
\end{thm}
In the limit case $ k \rightarrow 0$, we obtain $\capa(A) \geq \frac{l^2}{\area(A)}$.\\

If $C$ has constant curvature $K$, then the theorem says the following. Among all annuli of fixed area on a cylinder with baseline $\gamma$, the annulus with constant length, centered around $\gamma$ has minimal capacity (see \textit{Fig.~\ref{fig:negcurv2}}). In this case the inequality is sharp. This is also the annulus, whose boundary line has minimal length among all annuli of fixed area. This means that for this annulus the isoperimetric inequality for (non-contractible) annuli on cylinders is sharp. There is no such lower bound depending on the area for contractible annuli.
In the case of constant positive curvature, we obtain:
\begin{thm}
Let $C$ be a cylinder of constant curvature $K=k^2$ and let $\gamma$ be its baseline of length $l(\gamma)=l$. Let $A \subset C$ be a non-contractible annulus of finite area $\area(A)$. Then for $W'=\arcsin(\frac{k\area(A)-l\sin(ka)}{l})$ and $h_2(x)= \log\left(\frac{1+\sin(x)}{\cos(x)}\right)$, we have
\[
     \capa(A) \geq \frac{k \cdot l}{h_2(W') - h_2(-ka)}.
\]
\label{thm:main2}
\end{thm}
The theorem says that among all annuli of fixed area on a cylinder of constant positive curvature with baseline $\gamma$, the annulus with constant length, where one boundary is the shortest boundary of the cylinder itself, has minimal capacity (see \textit{Fig.~\ref{fig:poscurv2}}). As in the case of constant negative curvature this is the annulus that satisfies the isoperimetric inequality for annuli. For more general results about the isoperimetric inequality and capacity, see \cite{gr}.\\
The result for the cylinders of constant curvature is obtained with the help of a symmetrization process that results in an annulus of minimal capacity. The result for cylinders of variable negative curvature is obtained by comparing the capacity of an annulus with the capacity of an annulus on a comparison cylinder of higher, constant, curvature. We will prove these theorems in section 3. Using the same methods, we will derive lower and upper bounds of the capacities of annuli on cylinders in section 4. These results are stated in \textbf{Theorem \ref{thm:constant_k}}. The methods developed here have been applied to estimate the energy of harmonic forms on Riemann surfaces in \cite{mu}.\\
In \cite{g}, Gehring provides elementary estimates for the capacity of rings in $\R^n$, taken with respect to an arbitrary metric. Here the upper bound is constructed with the help of a test function that increases linearly along geodesics that realize the distance between a point on one boundary and the second boundary. It is noteworthy that this method can be adapted to obtain an upper bound for annuli on cylinders.\\

\section{Preliminaries}

\textit{Basic definitions}\\
\\
There exist two types of annuli on $C$, contractible and non-contractible annuli (see \textit{Fig.~\ref{fig:poscurv2}} and \textit{ Fig.~\ref{fig:negcurv2}}). Here a \textit{non-contractible  annulus} or shortly \textit{annulus} $A \subset C$ on a cylinder $C$ is a set that can be obtained by a continuous deformation of $C$. More precisely there exists an isotopy
\[
J:C\times[0,1] \rightarrow C, \text{ \ \ such that \ \ }  J(\cdot,0)=id  \text{ \ \  and  \ \ } J(C,1)=A.
\]
We denote by $\partial_1 A$ and $\partial_2 A$ the two connected boundary components that constitute the boundary $\partial A$. Let
\[
\phi: [0,c_1] \times ]c_2,c_3[ \mod \{(0,s) \sim (c_1,s) \mid s \in ]c_2,c_3[\}  \rightarrow C
\]
be a bijective parametrization of $C$ and $G$ the corresponding metric tensor. Let $F \in \lip(\bar{A})$ a Lipschitz function on the closure of $A$. Then the \textit{energy of $F$ on $A$}, $E_A(F)$ is given by
\[
    E_A(F) = \iint\limits_{\phi^{-1}(A) } {\left\| {D(F \circ \phi)} \right\|_{G^{-1}}^2 \sqrt{\det(G)}}.
\]
The \textit{capacity} of an annulus $A$, $\capa(A)$ is given by
\[
   \capa(A)= \mathop {\inf }\{E_A(F) \mid \{ F \in \lip(\bar{A}) \mid F |_{\partial_1 (A)} = 0 , F |_{\partial_2 (A)} = 1 \} \}.
\]
For further information about the definition of the capacity in metric spaces, see \cite{tr}.\\
\\
\textit{Fermi coordinates}\\
\\
Let $E_K$ be a simply connected surface of either constant curvature or variable negative curvature.
Let $\nu \subset E_K$ be a geodesic. $\nu$ divides $E_K$ into two parts, $E_K^+$ and $E_K^-$. Let $\gamma \subset \nu$ be  a geodesic arc of length $l(\gamma) < \diam(E_K)$ with endpoints $p_1$ and $p_2$.\\
It follows from \cite{kl}, p. 62-64 that for each point $p \in \gamma$ there exists a unique geodesic $\delta_p$ that is perpendicular to $\gamma$ and that passes $p$. The \textit{Fermi coordinates with base point $p_1$ and baseline $\gamma$} in $E_K$ are an injective parametrization
\[
\psi:D=[0,l(\gamma)] \times ]a,b[ \rightarrow E_K, \psi: (t,s) \mapsto \psi(t,s)
\]
with $a < 0 \leq b$ and $|b| \leq |a| \leq \frac{\diam(E_K)}{2}$, such that the parametrization satisfies the following conditions :\\
Each point $q=\psi(t,s) \in \psi(D)$ can be reached in the following way. Starting from the base point $p_1$ we first we move along $\gamma$ the directed distance $t$ to $\psi(t,0)=p$ and then from $\psi(t,0)=p$, we now move along $\delta_p$ the directed distance $s$ to $\psi(t,s)=q$. \\
An image $\psi(D) = S$ is a strip. We remind that a cylinder with baseline $\gamma$ is
\[
C= S \mod M,
\]
such that the side of $S$ containing $p_1$ is identified by an isometry $M \in \iso^+(E_K)$ with the side of the strip $S$ containing $p_2$.\\
Any cylinder of variable negative curvature or constant curvature can be parametrized in Fermi coordinates (see \cite{kl}, p. 62-64). It is a well-known fact that in the case of negative curvature the image of $\gamma$ is the shortest simple closed geodesic on $C$.\\
\\
\textit{Cylinders of constant curvature} \\
\\
If $E_K$ has constant curvature $K$ then a model of $E_K$ is the hyperbolic plane if $K < 0$, the Euclidean plane if $K=0$, and the sphere if $K>0$. We will use the following models for the simply connected surfaces of constant negative curvature and positive curvature : \\
In the case of negative curvature $K =-k^2$, set $E_K=\Hy$, where $\Hy$ is the Poincare model of the hyperbolic plane. It is the following subset of the complex plane $\C$ :
\[
   \Hy = \{ z = x + i y \in \C \mid y > 0 \}
\]
with the hyperbolic metric
\[
    ds^2 = \frac{1}{(ky)^2}(dx^2 + dy^2).
\]
As the metric is conformal, the hyperbolic energy of a function $F$ on a set $L$, $F: L \subset \Hy \rightarrow \R$ is equal to the Euclidean energy of $F$.\\
The parametrization $\psi : [0,l]\times \R \mod \{(0,s) \sim (l,s) \mid s \in \R \}  \rightarrow \Hy$
\[
   \psi(t,s) := \frac{\exp(kt)}{\cosh(ks)}(\sinh(ks) + i)
\]
parametrizes a cylinder $C$ with baseline $\gamma=\{i y  \mid y \in [1,\exp(kl)] \}$.\\
\\
In the case of positive curvature $K = k^2$, set $E_K=\Sp_k^2$, where $\Sp_k^2$ is the sphere $\partial B_{\frac{1}{k}}(0) \subset \R^3$, the boundary of the ball of radius $\frac{1}{k}$ . The parametrization $\psi : [0,l]\times ]a,b[ \mod \{(0,s) \sim (l,s) \mid s \in ]a,b[ \}  \rightarrow \Sp_k^2$
\[
   \psi(t,s) := \frac{1}{k}(\cos(kt)\cos(ks),\sin(kt)\cos(ks),\sin(ks))
\]
parametrizes a cylinder $C \subset \Sp_k^2$ with baseline $\gamma=\{(\frac{1}{k}\cos(ky),0,0)  \mid y \in [0,l] \}$.\\

Using these models, we obtain the metric tensor $G$ with respect to the Fermi coordinates :
\begin{equation}
G(t,s)=\left( {\begin{array}{*{20}c}
   {h(t,s)^2 } & 0  \\
   0 & 1  \\
\end{array} } \right),
\label{eq:tensorG}
\end{equation}
where
\begin{equation}
h(t,s) = h(s) = \cosh(ks), \text{ \ \ if \ \ } K=-k^2 \text{ \ \ and \ \ } h(t,s) = h(s) = \cos(ks), \text{ \ \ if  \ \ }K=k^2.
\label{eq:hs}
\end{equation}

\textit{Cylinders of variable negative curvature} \\
\\
Let $C=S \mod M$ be a cylinder with variable negative curvature bounded from above by $K$. Then with respect to the Fermi coordinates  $\psi$ for $C$, the metric tensor $G$ has the same shape as in equation (\ref{eq:tensorG}).\\
Let $\gamma$ be the baseline of $C$. Now the Fermi coordinates parametrize $\gamma$ with unit speed. It follows furthermore from the fact that $\gamma$ is a geodesic in $C$ that the derivative of $h(t,s)$ with respect to $s$ is zero, if $s$ equals zero. Hence
\[
h(t,0)=1 \text{ \ \ and \ \ }  \frac{\partial h(t,0)}{\partial s} = 0.
\]
If $K(t,s)$ is the curvature in the point $\psi(t,s) \in C$ then it follows from equation (\ref{eq:tensorG}) that
\[
   h(t,s)\cdot K(t,s)= -\frac{\partial^2 h(t,s)}{\partial s^2}.
\]
If the curvature of $C$ is bounded between $-k_2^2$ and $-k_1^2 =K$  then we obtain the following result with the help of the last two equations. Applying Cauchy's mean value theorem twice to $\frac{h(t,s)-h(t,0)}{\cosh(k_i s) -\cosh(k_i 0)}$ for $i \in \{1,2\}$, we obtain :
\begin{equation}
\cosh(k_2s) \geq h(t,s) \geq \cosh(k_1 s) \text{ \ for all \ } t\in [0,l(\gamma)].
\label{eq:ineqk1k2}
\end{equation}
Using this inequality and the previous two equations a simple application of the mean value theorem of differentiation applied to $\frac{\frac{\partial h(t,s)}{\partial s}- \frac{\partial h(t,0)}{\partial s}}{s-0}$  gives us for all $t\in [0,l(\gamma)]$
\[
\frac{\partial h(t,s)}{\partial s} < 0 \text{ \ \ for \ \ } s < 0  \text{ \ \ and \ \ }   \frac{\partial h(t,s)}{\partial s} > 0 \text{ \ \ for \ \ } s > 0.
\]

\textit{Area and  directional capacity}\\
\\
Using the formulas for the metric tensor $G$, the area of an annulus $A$, $\area(A)$ is given by
\begin{equation}
\area(A)=\iint\limits_{\psi^{-1}(A)} {\sqrt{\det(G(t,s))}\,ds \, dt }=\iint\limits_{\psi^{-1}(A)} {h(t,s)  \,ds \, dt }.
\label{eq:area}
\end{equation}

Now let $A \subset C$ be an annulus and $F:A \rightarrow \R$ be a $\lip(\bar{A})$ function. For a $x=\psi(t_0,s_0) \in C$ denote by
\[
p_{\psi} : T_x C \rightarrow \{\lambda \cdot \frac{\partial \psi(t_0,s_0)}{\partial s} \mid \lambda \in \R \}
\]
the orthogonal projection of a tangent vector in $x$ onto the subspace spanned by $\frac{\partial \psi(t_0,s_0)}{\partial s}$. We denote by $E_A(\partial_2 F)= E_A(p_{\psi}(DF))$ the energy of this orthogonal projection of $DF$.\\
Let $L \subset C$ be a strip or annulus. For technical purposes, we also define the \textit{capacity of $L$ in direction $\partial_2$}, $\caps(L)$ or shortly \textit{directional capacity} by
\[
   \caps(L)= \mathop {\inf }\{E_L(\partial_2 F) \mid \{ F \in \lip(\bar{L}) \mid F |_{\partial_1 (L)} = 0 , F |_{\partial_2 (L)} = 1 \}. \}
\]
Using Fermi coordinates, we obtain for the energy $E_A(F)$ of $F$ on $A$, with $F \circ \psi =f$ :
\begin{equation}
E_A(F)=\iint\limits_{\psi^{-1}(A)}{\frac{\frac{\partial f(t,s)}{\partial t}^2}{h(t,s)}+h(t,s)\frac{\partial f(t,s)}{\partial s}^2 \,ds \, dt } \geq \iint\limits_{\psi^{-1}(A)} {h(t,s)\left(\frac{\partial f(t,s)}{\partial s}\right)^2  \,ds \, dt } = E_A(\partial_2 F).
\label{eq:main}
\end{equation}
Using equation (\ref{eq:main}), we obtain the following lemma. Note that we can drop the condition $a < 0$.
\begin{lem} Let $S=\psi([0,l(\gamma)]\times ]a,b[) \subset E_{K}$ be a strip and $C = S \mod M$ be a cylinder with baseline $\gamma$. Then we have
\[
\caps(C)=\caps(S) = \int \limits_{t=0}^{l(\gamma)}{  \frac{1}{H(t,b)-H(t,a)} \,dt} , \text { \ \ where \ \ }
H(t,s)= \int {\frac{1}{h(t,s)} \,ds }.
\]
\label{thm:tech}
\end{lem}

\textbf{proof of Lemma~\ref{thm:tech}}  For any function $F \in \lip(\bar{C}), F|_{\partial_1 C} =0, F|_ {\partial_2 C}=1$, with $F \circ \psi =f$, we obtain by inequality (\ref{eq:main}) that
\[
E_{C}(F)\geq \iint\limits_{\psi^{-1}(C)} {h(t,s)\left(\frac{\partial f(t,s)}{\partial s}\right)^2  \,ds \, dt }.
\]
Solving the Euler-Lagrange equation (see \cite{ge} p. 152-154), we can determine the function $P=p \circ \psi^{-1} \in \lip(\bar{C})$ that satisfies the boundary conditions on $C$ and such that $p$ minimizes the second integral in the above inequality (\ref{eq:main}).
We obtain
\begin{eqnarray*}
p(t,s)&=&c_1 H(t,s)-c_2 , \text{ \ \ where \ \ } \\
c_1&=&\frac{1}{H(t,b)-H(t,a)} \text{ \ \ and  \ \ } c_2=\frac{H(t,a)}{H(t,b)-H(t,a)}.
\end{eqnarray*}

Hence for all $F \in \lip(\bar{C}), F|_{\partial_1 C} =0, F|_{\partial_2 C}=1$ :
\[
E_{C}(F)\geq \iint\limits_{\psi^{-1}(C)} {h(t,s)\left(\frac{\partial f(t,s)}{\partial s}\right)^2  \,ds \, dt } = \int\limits_{t=0}^l {\int\limits_{s=a}^b {h(t,s)\left(\frac{\partial p(t,s)}{\partial s}\right)^2  \,ds }\, dt }.
\]

We obtain $\caps(C) = \caps(S)= \int \limits_{t=0}^{l}{  \frac{1}{H(t,b)-H(t,a)} \,dt} $, from which follows the lemma. $\square$\\

If $E_K$ has constant curvature, then we obtain :
\begin{equation}
H(t,s) = H(s) = \left\{ {\begin{array}{*{20}c}
   \frac{h_1(ks)}{k}  \\
   \frac{h_2(ks)}{k}  \\
\end{array} } \right. \text{ if }
\begin{array}{*{20}c}
    K = -k^2 \\
   K = k^2  \\
\end{array}.
\label{eq:Hconst}
\end{equation}
Here $h_1(s)=2\arctan(\exp(s))$ and $h_2(s)=\log\left(\frac{1+\sin(s)}{\cos(s)}\right)$ (see \textbf{Theorem~\ref{thm:main1}} and \textbf{\ref{thm:main2}}).\\
In the case of constant curvature $P$ is the minimizing function for the capacity problem and we have
\begin{equation}
\capa(C)=\caps(S)=E_S(P), \text{ \ where \ } P \text{ \ harmonic and \ } \frac{\partial P\circ \psi}{\partial t}=0.
\label{eq:capconst}
\end{equation}
We note that in any case $\caps(S)$ is decreasing, if the length of the strip $S$ is increasing.

\section{Lower bounds on the capacity depending on the area}

\textit{Cylinders of constant curvature}\\
\\
We will first prove \textbf{Theorem \ref{thm:main1}} and \textbf{\ref{thm:main2}} for cylinders of constant curvature. The proof consists in replacing in three steps the initial annulus $A=A_0$, with a set $A_i$ such that
\[
 \caps(A_i) \geq  \caps(A_{i+1}) \text{ \ \ and \ \ } \area(A_i)=\area(A).
\]
In the final step, we obtain an annulus $A_3$ of constant length that satisfies $\caps(A_{3})=\capa(A_{3})$.\\

\textbf{Step 1 : Reduction of the number of subsections}\\

Let $A$ be an annulus of fixed area on a cylinder $C$ with baseline $\gamma$. For each $p \in \gamma$ there is a geodesic arc $\delta_p$ that is perpendicular to $\gamma$ and that passes through $p$. We call  $\delta_p \cap A =\eta_p$ a \textit{section} of $A$. We call a  geodesic arc $\eta_i$ that forms a connected component of a section a \textit{subsection}. We furthermore define the type of the subsection. The type depends on the boundary condition imposed by the capacity problem. We say that a subsections \textit{$\eta_i$ is of type $aa$}, if the boundary conditions on $\eta_i$ imply that the boundary values on both sides of the subsection are equal. We say that \textit{$\eta_i$ is of type $ab$} if the boundary conditions imply that the values on the two sides are different. Note that $\caps(\cdot)$ is well-defined on sections and subsections.\\
We may assume that a section $\eta_p$ is given in Fermi coordinates by
\[
\eta_p = \bigcup\limits_{i=1}^{n_p} \eta_i = \psi(\{t\}  \times \bigcup \limits_{i=1}^{n_p} ]a_i,b_i[ )  \text{ \ \, where \ \ } a_i < b_i < a_{i+1}.
\]
Here the number of subsections is finite. It will be clear from the proof that the proof also applies if the number is infinite.\\
In this step we reduce the number of subsections in each section of $A$ to one, such that the total area of $A$ does not change. Note that in each section there must be a subsection of type $ab$, because $A$ is a non-contractible annulus.\\
To execute \textbf{Step 1} we apply the following algorithm to each section $\eta_p = \bigcup \limits_{i=1}^{n_p} \eta_i $  :\\
Starting from $i=1$, we go consecutively through the subsections $\eta_i$. For fixed $i$ we first replace $\eta_i \cup \eta_{i+1}$ by one subsection $\eta'$. We obtain $\eta'$ by elongating one of the subsection in direction of the other. If both $\eta_i$ and $\eta_{i+1}$ are of the same type, we elongate $\eta_i$ in direction of $\eta_{i+1}$.\\
If one subsection is of type $aa$ and the other subsection is of type $ab$, then we elongate the subsection of type $ab$ in direction of the subsection of type $aa$.\\
Let WLOG $\eta_i$ be the subsection that is being elongated. We note that for a subsection $\eta_i$ of type $aa$ we have $\caps(\eta_i)=0$. It follows from \textbf{Lemma~\ref{thm:tech}} that $\caps(\eta_i)$ decreases, if the length of $\eta_i$ increases. Hence
\[
\caps(\eta_i \cup \eta_{i+1} )  = \caps(\eta_i) + \caps(\eta_{i+1})  \geq \caps(\eta_i) \geq \caps(\eta').
\]
We elongate $\eta_i$ to obtain $\eta'=\psi(\{t\}\times ]a_i,b[)$ in a way such that
\begin{equation}
\int\limits_{a_i}^{b_i} {h(t,s) \, ds } +  \int\limits_{a_{i+1}}^{b_{i+1}} {h(t,s) \, ds } =  \int\limits_{a_i}^{b} {h(t,s) \, ds } .
\label{eq:redu_elong}
\end{equation}
Using (\ref{eq:redu_elong}) in equation (\ref{eq:area}), this condition ensures that the set $A_1$ obtained in this way has the same area as the annulus $A$.\\
Finally we set $\eta':=\eta_{i+1}$ and go to ${i+1}$.\\
When we apply this algorithm to a section $\eta_p$ we obtain a single connected section $\eta'_{p}$. We obtain $A_1$ by applying this algorithm to all sections of $A$. We have
\[
\caps(A) \geq \caps(A_1)  \text{ \ \ and \ \ } \area(A_1) = \area(A).
\]
\textbf{Step 2 : Positioning the sections}\\

A section $\eta_p$ of $A_1$ in Fermi coordinates is WLOG given by
\[
\eta_p = \psi(\{t\} \times ]a_1(t),a_2(t)[) .
\]
In this step we replace each $\eta_p$ with a section $\eta'_p$, such that the resulting annulus $A_2$ has the same area.\\
If $K>0$, we replace $\eta_p$ with a section $\eta'_p=\psi(\{t\} \times ]a,W(t)[)$ situated near the shorter boundary of $C$. If $K<0$, we replace it with a section $\eta'_p=\psi(\{t\} \times ]-W(t),W(t)[)$ centered around $\gamma$.\\
To ensure that the resulting annulus $A_2$ has the same area as $A_1$, we have to solve the equations
\begin{equation}
\int\limits_{a_1(t)}^{a_2(t)} {h(s) \, ds }  =  \int\limits_{-W(t)}^{W(t)} {h(s) \, ds }  \text{ \ \ if \ \ } K < 0 \text{ \ \ \ and \ \ \ }
\int\limits_{a_1(t)}^{a_2(t)} {h(s) \, ds }  =  \int\limits_{a}^{W(t)} {h(s) \, ds } \text{ \ \ if \ \ } K > 0
\label{eq:area_comp1}
\end{equation}
with respect to $W(t)$ to obtain the length coordinates of the section $\eta'_p$ (see equation (\ref{eq:hs}) and (\ref{eq:area})). Then it follows that
\begin{equation*}
\caps(\eta_p)  \geq \caps(\eta'_p)
\label{eq:step2_const}
\end{equation*}
from \textbf{Lemma~\ref{thm:tech}} and equation (\ref{eq:Hconst}).\\
Replacing all sections $\eta_p$ of $A_1$ by the corresponding section $\eta'_p$ we obtain the annulus $A_2$ of lower directional capacity with $\area(A_2)=\area(A)$.\\
\\
\textbf{Step 3 : Averaging the lengths of the sections}\\
\\
In this step we replace $A_2$ by an annulus $A_3$ of constant length. If $K<0$ we replace $A_2$ with an annulus centered around $\gamma$. If $K>0$ we replace $A_2$ with an annulus situated near the shorter boundary of $C$.\\
A section of $A_2$ in Fermi coordinates is WLOG given by
\[
\psi([\{t\} \times ]a_1(t),a_2(t)[), \text{ \ where \ }  -a_1(t)=a_2(t)=W(t), \text { \ if \  } K<0 \text { \ and \ } a_1(t) = a, \text{ \ if \ } K >0 .
\]
To obtain an annulus $A_3$ of the same area as $A_2$ we have to solve the equation
\begin{align*}
\area(A_2)&= \int \limits_{t=0}^l \int \limits_{s=-W(t)}^{W(t)} {{h(s) \, ds } \, dt} &=  \int \limits_{t=0}^l \int\limits_{s=-W}^{W} {{h(s) \, ds } \, dt } &= \area(A_3)  &\text{ \ \ if \ \ } K < 0,\\
\area(A_2)&= \int \limits_{t=0}^l \int \limits_{s=a}^{W(t)} {{h(s) \, ds } \, dt} &=  \int \limits_{t=0}^l  \int \limits_{s=a}^{W} {{h(s) \, ds } \, dt } &= \area(A_3) &\text{ \ \ if \ \ } K > 0.
\end{align*}
with respect to $W$ to obtain the length coordinates of $A_3$. We obtain :
\begin{align}
\nonumber
W&=\frac{1}{k}\arcsinh\left(\frac{k\area(A)}{2 \cdot l}\right)  &\text{ \ \  if \ \ } K < 0,\\
\label{eq:W3}
W&=\frac{1}{k}\arcsin\left(\frac{k\area(A)+l\sin(ka)}{l}\right) &\text{ \ \ if \ \ } K > 0.
\end{align}
It remains to show that
\[
   \caps(A_3) \leq \caps(A_2).
\]
Applying \textbf{Lemma~\ref{thm:tech}} and equation (\ref{eq:Hconst}), this is equal to
\begin{align}
\nonumber
\frac{k \cdot l}{h_1(kW)-h_1(-kW)} &\leq \int \limits_{t=0}^l { \frac{k}{h_1(kW(t))-h_1(-kW(t))}  \, dt }   &\text{ \ \ if \ \ } K < 0,\\
\label{eq:caps3}
\frac{k \cdot l}{h_2(kW) - h_2(ka)} &\leq \int \limits_{t=0}^l  { \frac{k}{h_2(kW(t)) - h_2(ka)} \, dt}  &\text{ \ \ if \ \ } K > 0.
\end{align}

This result follows from an application of the integral version of Jensen's inequality. If $\varphi : \R_+ \rightarrow \R$ is a convex function, then :
\begin{equation*}
\varphi\left( \, \, \int \limits_{t=0}^{l}{ f(t) \, dt}\right) \leq   \int \limits_{t=0}^l { \frac{1}{l} \varphi(l \cdot f(t)) \, dt}.
\end{equation*}

Here we use the result for $W$  from equation (\ref{eq:W3}) in inequality (\ref{eq:caps3}) and set
\begin{align*}
f(t)&:= \int\limits_{s=-W(t)}^{W(t)} h(s) \, ds \text{ \ \ and \ \ } &\varphi(x)& := \frac{k\cdot l}{h_1(\arcsinh(\frac{k \cdot x }{2 \cdot l}))-h_1(-\arcsinh(\frac{k \cdot x }{2 \cdot l}))} &\text{ \ \ if \ \ } K < 0,\\
f(t)&:= \int\limits_{s=a}^{W(t)} h(s) \, ds \text{ \ \ and \ \ } &\varphi(x)& := \frac{k \cdot l}{h_2(\arcsin(\frac{k\cdot x+l\sin(ka)}{l}) ) - h_2(ka)} &\text{ \ \ if \ \ } K > 0
\end{align*}
in Jensen's inequality. This way we obtain inequality (\ref{eq:caps3}).\\
We obtain always the same annulus $A_3$, independent of the starting annulus $A$. This annulus satisfies $\caps(A_3)=\capa(A_3)$ (see (\ref{eq:capconst})). This proves \textbf{Theorem \ref{thm:main1}} and \textbf{\ref{thm:main2}} in the case of cylinders of constant curvature.  \ \ \ \ \  $\square$\\
\\
\textit{Cylinders of variable negative curvature}\\
\\
Now let $C$ be a cylinder of variable negative curvature with baseline $\gamma$ of length $l$. We assume that the curvature of $C$ is bounded from above by $K=-k_1^2$. Let $A=A_0$ an annulus of fixed area. Here we proceed in two steps.\\
First we apply \textbf{Step 1} from the case of constant curvature to obtain $A_1$. From \textbf{Step 1} we obtain an annulus $A_1$ of equal area and with lower directional capacity than $A_0$.\\
In the following \textbf{Step 2} we show that $\caps(A_1)$ is always bigger than the capacity of a certain comparison annulus of fixed length on a comparison cylinder of constant curvature $-k_1^2$.\\
\\
\textbf{Step 2 : Comparison with an annulus of constant negative curvature}\\
\\
Each section of $A_1$ consists of a single geodesic arc $\eta_p$. Now we compare the directional energy $\caps(\eta_p)$ with the directional energy of a section $\eta'_p$ in a comparison cylinder $C'$ with baseline $\gamma$ of length $l(\gamma)=l$ of constant negative curvature $-k_1^2$.\\
We have that $\eta_p$ is given in Fermi coordinates by
\[
\eta_p = \psi(\{t\} \times ]a_1(t),a_2(t)[) \subset A_1.
\]
We choose a section $\eta'_p \subset C'$ centered around $\gamma$, such that for $t \in [0,l]$
\[
    \int \limits_{s=a_1(t)}^{a_2(t)} { h(t,s)} \,ds = I(t) =  \int \limits_{s=-W(t)}^{W(t)} { \cosh(k_1 s) } \,ds.
\]
This condition ensures that the corresponding annulus $A_2 \subset C'$ has the same area as $A_1$ (see equation (\ref{eq:hs}) and (\ref{eq:area})).
The above equation is equal to
\begin{equation}
I(t) = \frac{2}{k_1}\sinh(k_1 \cdot W(t)) \Leftrightarrow W(t)=\frac{1}{k_1}\arcsinh(\frac{k_1 \cdot I(t)}{2}).
\label{eq:It}
\end{equation}

To obtain $\caps(A_1) \geq \caps(A_2)$ with the help of \textbf{Lemma \ref{thm:tech}} we have to show that for all $t \in [0,l]$
\begin{equation}
\frac{1}{\int\limits_{a_1(t)}^{a_2(t)} {\frac{1}{h(t,s)} \,ds}} \geq \frac{1}{\int\limits_{-W(t)}^{W(t)} {\frac{1}{\cosh(k_1 s)} \,ds}} \Leftrightarrow  \int\limits_{a_1(t)}^{a_2(t)} {\frac{1}{h(s,t)} \,ds} \leq \int\limits_{-W(t)}^{W(t)} {\frac{1}{\cosh(k_1 s)} \,ds}.
\label{eq:inequ_varcurv}
\end{equation}
We now replace $W(t)$ from equation (\ref{eq:It}) in equation (\ref{eq:inequ_varcurv}). We note that $-\arcsinh(x)=\arcsinh(-x)$. Using integration by substitution we obtain :
\[
\int\limits_{-W(t)}^{W(t)} {\frac{1}{\cosh(k_1 s)} \,ds} =\frac{1}{2} \int\limits_{-I(t)}^{I(t)} {\frac{1}{\cosh^2(\arcsinh(\frac{k_1 \cdot s}{2}))} \,ds}.
\]
with $\cosh^2(\arcsinh(x)) = 1+x^2$ this simplifies to
\[
 \int\limits_{-W(t)}^{W(t)} {\frac{1}{\cosh(k_1 s)} \,ds} =\frac{1}{2} \int\limits_{-I(t)}^{I(t)} {\frac{1}{1+(\frac{k_1 \cdot s}{2})^2} \,ds}.
\]
Now as $\frac{1}{a}\arctan(a x)'=\frac{1}{1+(ax)^2}$ and as $-\frac{1}{a}\arctan(a x)= \frac{1}{a}\arctan(a x)$ we obtain that
\[
 \int\limits_{-W(t)}^{W(t)} {\frac{1}{\cosh(k_1 s)} \,ds} = \frac{2}{k_1}\arctan(\frac{k_1 I(t)}{2}).
\]
By inequality (\ref{eq:ineqk1k2}) we have that $h(t,s) > \cosh(k_1 s)$ for all $t \in [0,l]$. Hence we obtain on the left-hand side of (\ref{eq:inequ_varcurv})
\[
\int\limits_{a_1(t)}^{a_2(t)} {\frac{1}{h(t,s)} \,ds} \leq  \int\limits_{a_1(t)}^{a_2(t)} {\frac{1}{\cosh (k_1 s)} \,ds} = \frac{2}{k_1}\left(\arctan(\exp(k_1 a_2(t))) - \arctan(\exp(k_1 a_1(t)))\right).
\]
We now consider the right-hand side of inequality (\ref{eq:inequ_varcurv}). Due to the monotonicity of the $\arctan$ function and as $h(t,s) > \cosh(k_1 s)$ we have
\[
\frac{2}{k_1}\arctan(\frac{k_1 I(t)}{2}) \geq \frac{2}{k_1}\arctan\left(\frac{k_1}{2} \int\limits_{a_1(t)}^{a_2(t)} \cosh (k_1 s) \,ds\right) = \frac{2}{k_1}\arctan\left(\frac{ \sinh(k_1 a_2(t))-\sinh(k_1 a_1(t))}{2}\right).
\]
As for all $a_1(t),a_2(t) \in \R, a_2(t) >a_1(t)$
\begin{equation*}
\arctan(e^{(k_1 a_2(t)}) - \arctan(e^{k_1 a_1(t)})  \leq \arctan\left(\frac{ \sinh(k_1 a_2(t))-\sinh(k_1 a_1(t))}{2}\right) ,
\label{eq:step2_var}
\end{equation*}
we obtain in total inequality (\ref{eq:inequ_varcurv})
\[
\int\limits_{a_1(t)}^{a_2(t)} {\frac{1}{h(t,s)} \,ds}  \leq \int\limits_{-W(t)}^{W(t)} {\frac{1}{\cosh(k_1 s)} \,ds}.
\]
It follows that $\caps(A_1) \geq \caps(A_2)$. Hence for all $A_1$ from \textbf{Step 1} there exists an comparison annulus $A_2$ in a cylinder of constant curvature $-k_1^2$, such that $\area(A)=\area(A_1)=\area(A_2)$ and such that $A_2$ has lower directional capacity. For this annulus we obtain a lower bound based on the area from the previous proof. This concludes the proof in the case of annuli on cylinders of variable negative curvature. $\square$

\section{Inequalities for the capacity of an annulus on a cylinder}

The upper and lower bound for the capacity of an annulus $A$ can be obtained in the following way. To obtain an upper bound, we can evaluate the energy of any test function $F_T \in \lip(\bar{A})$ that satisfies the boundary conditions for the capacity problem for $A$. We can easily construct a test function by adjusting the minimizing function $P$ from \textbf{Lemma \ref{thm:tech}} to the boundary, such that $P$ is the minimizing function for $\caps(A)$. We can evaluate our choice by evaluating $\caps(A)$, which provides the lower bound for $\capa(A)$. This approach works immediately for an annulus $A$ that can be parametrized in Fermi coordinates in the following way :
\[
A=\psi\{ (t,s) \mid s \in [a_1(t),a_2(t)],t \in [0,l] \},
\]
where $a_1(\cdot)$ and $a_2(\cdot)$ are piecewise derivable functions with respect to $t$.\\
In this case, we say that our annulus is of \textit{type A}. As $a_1(\cdot)$ and $a_2(\cdot)$ are functions each section of $A$ consists of a single geodesic arc.\\
We say an annulus is of \textit{type B} if it is not of type A and if its boundary is a piecewise differentiable curve. In this case the approach can be adapted to obtain an upper or lower bound. Here the lower bound can be constructed by the same method, which we present in the following. Though the method can also be adapted to obtain an upper bound for any annulus of type B, we think that this upper bound deviates too much from the real value of the capacity and we will not present this approach here.\\
As mentioned in the introduction, different test functions are used in \cite{g} to obtain an upper bound for the capacity of annuli. These methods can also be applied in the current situation. However, the upper bound presented in this section is more practical and explicit. The result of this section has been used to estimated the energy of harmonic forms dual to a canonical homology basis on Riemann surfaces in \cite{mu}.\\
\\
\textit{Annuli of type A}\\
\\
The following theorem can be concluded from the discussion above. For a definition of $H$ see \textbf{Lemma~\ref{thm:tech}} :\\

\begin{thm}
Let $S \subset E_K$ be a strip and $C = S \mod M$ be a cylinder with baseline $\gamma$ of length $l(\gamma)=l$. Let $A  \subset C$ be an annulus of type A and $P \in \lip(\bar{A})$ be the function whose energy realizes $\caps(A)$, then
\[
       E_A(P) \geq  \capa(A) \geq \caps(A) = E_A(\partial_2P) = \int\limits_{t=0}^{l} {\frac{1}{H(t,a_2(t))-H(t,a_1(t))} \, dt}.
\]
If $C$ has constant curvature $K$ then we obtain with $q_i(t)=\frac{\partial H(s_0)}{\partial s}|_{s_0=a_i(t)}\cdot a_i'(t)$ for $i \in \{1,2\}$:
\[
\int\limits_{t=0}^{l} {\frac{1+\frac{q_1(t)^2+q_1(t)q_2(t)+q_2(t)^2}{3}}{H(a_2(t))-H(a_1(t))} \, dt } \geq \capa(A) \geq \int\limits_{t=0}^{l} {\frac{1}{H(a_2(t))-H(a_1(t))} \, dt}.
\]
\label{thm:constant_k}
\end{thm}

\textbf{proof of Theorem \ref{thm:constant_k}}  The first inequality follows from the previous paragraph. The lower bound follows from the representation of the optimal function $P$ for $\caps(A)$ in Fermi coordinates and \textbf{Lemma~\ref{thm:tech}}.\\
It remains to prove the first part of the second inequality. For the upper bound, we will calculate $p=P\circ \psi^{-1}$ explicitly.  For fixed $t \in [0,l]$, the boundary conditions imply
\begin{eqnarray*}
p(t,s)&=&c_1(t) H(s)-c_2(t) , \text{ \ \ where \ \ } \\
c_1(t)&=&\frac{1}{H(a_2(t))-H(a_1(t))} \text{ \ \ and  \ \ } c_2(t)=\frac{H(a_1(t))}{H(a_2(t))-H(a_1(t))}.
\end{eqnarray*}

We now calculate the energy of $P$ on $A$, $E_A(P)$. It follows with $H(s)'=\frac{1}{h(s)}$ :

\[
\frac{\partial p(t,s)}{\partial s} = \frac{c_1(t)}{h(s)}   \text{ \ \ and \ \ }
\frac{\partial p(t,s)}{\partial t} = c_1'(t)\cdot H(s)-c_2'(t).
\]

Hence we obtain :

\[
   E_A(P)=\int\limits_{t=0}^{l} {\int\limits_{s=a_1(t)}^{a_2(t)} {{(c_1'(t)\cdot H(s)-c_2'(t) )}^2 \cdot \frac{\partial H(s)}{\partial s} + {c_1(t)}^2 \cdot \frac{\partial H(s)}{\partial s} \,dt} \,ds}.
\]

Evaluating the integral with respect to $s$, we have

\[
E_A(P)=\int\limits_{t = 0}^l {\left. {\frac{{(c_1 '(t)u - c_2 '(t))^3 }}
{{3c_1 '(t)}} + c_1 (t)^2 u} \right|_{u =  H(a_1(t))}^{H(a_2(t))} dt}.
\]

We have for $c_1'(t)$ and $c_2'(t)$, as $q_i(t)=\frac{\partial H(s_0)}{\partial s}|_{s_0=a_i(t)}\cdot a_i'(t)$ for $i \in \{1,2\}$:
\[
c_1'(t) = \frac{q_1(t) - q_2(t)}{(H(a_2(t))-H(a_1(t)))^2} \text{ \ \ and \ \ }
c_2'(t)= \frac{q_1(t)H(a_2(t)) - q_2(t)H(a_1(t))}{(H(a_2(t))-H(a_1(t)))^2}.
\]
With these equations $E_B(p)=E_A(P)$ simplifies to $\displaystyle{\int\limits_{t = 0}^l {\frac{1+\frac{q_1(t)^2 +q_1(t)q_2(t) + q_2(t)^2}{3}}{H(a_2(t))-H(a_1(t))} dt}.}$

This is the upper bound in \textbf{Theorem \ref{thm:constant_k}}. \ \ \ \ $\square$ \\

It is clear that the upper and lower bound are nearly optimal, if the boundary has only small variation, i.e. if $\displaystyle{ \int\limits_{t = 0}^l {|a_1'(t)|^2 + |a_2'(t)|^2 } \,dt }$ is small. If the variation is large, it might be possible to choose an annulus $A''$ in the interior of $A$, whose boundary line varies less. Then a test function can be constructed on this annulus as above. We can then evaluate the energy of this test function, to obtain a better upper bound.\\

\end{document}